# The Non-Substitution Theorem, Uniqueness of Solution and Convex combinations of basic optimal solutions for linear optimization

by

Somdeb Lahiri

ORCID: https://orcid.org/0000-0002-5247-3497

(Formerly) PD Energy University, Gandhinagar (EU-G), India

somdeb.lahiri@gmail.com

July 1, 2024.

This version: January 5, 2024.

## Abstract

Our first result is a statement of a somewhat general form of a non-substitution theorem for linear programming problems, along with a very easy proof of the same. Subsequently, we provide an easy proof of theorem 1 in Mangasarian (1979), based on a new result in terms of two statements that are each equivalent to a given solution of a linear programming problem being its unique solution. We also provide a simple proof of the result that states that the set of optimal solutions of a bounded linear optimization problem is the set of all convex combinations of its basic optimal solutions and the set of basic optimal solutions are the extreme points of the set of optimal solutions. We do so by appealing to Farkas' lemma and the well-known result that states that if a linear optimization problem has an optimal solution, it has at least one basic optimal solution. Both results we appeal to have easy proofs. We do not appeal to any version of the "Klein-Milman Theorem" or any result in advanced polyhedral combinatorics to obtain our results. As an application of this result, we obtain a simple proof of the Birkhoff-von Neumann Theorem.

**1. Introduction:** While much is known about the Non-Substitution Theorem and related issues in the case of activity analysis (see for instance Fujimoto, Silva and Villar (2003), Fujimoto, Herrrero, Ranade, Silva and Villar (2003), Lahiri (2023) and references in them), there seems to be some communication gap about what may be a "very general" form of the Non- Substitution Theorem for linear programming (LP) problems. This gap is vaguely reflected in the work by De Guili and Giorgi (2010) and some of the more significant works cited there. Our first result is a statement of a somewhat general form of such a non-substitution theorem along with a very easy proof of the same. It turns out that such a theorem does not rely on any "basic feasible solution" or "basic optimal solution" or for that matter on issues related to the non-degeneracy of any solution. The reason for this is, the simple fact that the Complementary Slackness Theorem of linear programming is concerned with the optimal solutions of the primal linear programming problem and its dual, and does not require such solutions to be basic.

Our next purpose in this paper, is to provide an easy proof of theorem 1 in Mangasarian (1979), without drifting into elaborate systems of inequalities. Theorem 1 of Mangasarian (1979), provides a necessary and sufficient condition for a given optimal solution for a LP problem, to be its unique solution. We prove this result using a proposition which to the best of our knowledge is an original result, though the motivation for it being theorem 1 in Mangasarian (1979), cannot be denied. The proposition is derived from a lemma that provides a necessary and sufficient condition for the zero vector to be the unique solution of the linear programming problem derived from the original one, by replacing the right-hand side of the linear equations in the constraints by the zero vector, and imposing the non-negativity constraints only on those variables which were zero in the given solution. The other variables are unconstrained in sign. The proposition states that the two equivalent statements in the lemma are equivalent to a given optimal solution of a LP problem being its unique optimal solution. Theorem 1 in Mangasarian (1979), which is the second theorem here, is an easy consequence of the proposition.

In section R4.4 of Lancaster (1968) we find the statements and proofs of two results, the second one being the finite dimensional version of the "Klein-Milman" theorem. Immediately preceding the statements of the two results the author provides his reasons for including them with the following words: "There are two theorems of considerable importance concerning extreme points. These theorems are fundamental to the theory of linear programming". The first sentence is unconditionally undisputable. We provide a simple proof of the result in LP, where the author may have intended the application of the two theorems mentioned above. This result, which is theorem 3 in this paper, states that the set of optimal solutions of a LP problem, is the set of all convex combinations of the set of all basic optimal solutions of the problem, and the set of all basic optimal solutions is the set of extreme points of the set of optimal solutions. Our proof relies on well-known results in LP that are available with extremely simple proofs in Lahiri (2020). A corollary of theorem 3 is that if the set of non-negative solutions of a system of linear equations is bounded, then it is the convex hull of its extreme points, and the extreme points are the non-negative basic solutions of the system of linear equations.

Simple proofs of theorem 3 and its corollary based on solving systems of linear equations, are available in Vajda (1961), as a combination of theorems 2.3 and 2.4 and a combination of theorem 2.1 and 2.2 in the book, respectively. The proof of our theorem 3 that is closest to any that we have encountered so far is the proof of theorem 3.9 in Goemans (2017), the latter being related but different from our theorem 3 in this paper. However, instead of relying on the fact that every solvable bounded linear optimization problem with a linear objective function has an optimal solution that is basic (Proposition 4 in Topic 2 of Lahiri (2020)) as we do here, to complete the proof of theorem 3.9 in Goemans (2017), the proof requires a Corollary 3.8 that is proved there. This corollary is a contribution to (somewhat advanced?) polyhedral combinatorics. Our proof of our theorem 3, is considerably simpler for those familiar with the comparatively elementary level results proved in the first three topics of Lahiri (2020). As an application of the corollary of theorem 3, we prove the Birkhoff-von Neumann theorem that every doubly stochastic matrix is a convex combination of permutation matrices. For the purpose, we use the proposition that any doubly stochastic matrix that is not a permutation matrix, can be expressed as a convex combination of two distinct doubly stochastic matrices, that are different from the one under consideration.

Sections 3, 4, 6 and sections 3, 4, 7 are each independent of one another. In section 5 we discuss results in section 4, and section 7 is an application of results obtained in section 6.

**2. The framework of analysis:** To begin with we shall formulate the LP problem in Mangasarian (1979) in its standard form as in Lucchetii, Radrizzani and Villa (2008).

Given any finite dimensional column vector z, we will represent its transpose by $z^T$ which is a row vector of the same dimension. Further given any two finite dimensional column vectors z, w, we will use $z \geq w$ to denote each coordinate of z is greater than or equal to the corresponding coordinate of w.

**Note:** In what follows we will implicitly assume that that the dimensions of matrices and vectors are conformable for multiplication, whenever such multiplication is required.

Thus, given an m×n matrix A, an n-dimensional column vector p and an m-dimensional column vector b, consider the LP problem:

Maximize $p^T x$, subject to $Ax = b$, $x \geq 0$.

*We shall refer to this problem* as **P1**.

If $x^*$ solves P1, then we will refer to $x^*$ as an **optimal solution** for/of P1.

For $j \in \{1, \ldots, n\}$, let $A^j$ denote the $j^{th}$ column of A.

**Observation:** In this context, an observation that is motivated by a very significant computational investigation reported in Serafini (2005) is the following.

Suppose $A^-$, $A^+$ are two m×n matrices, whose $(i, j)^{th}$ entry for $i \in \{1, \ldots, m\}$ and $j \in \{1, \ldots, n\}$ denoted by $a^-_{ij}$, $a^+_{ij}$ respectively satisfy $a^-_{ij} \leq a^+_{ij}$.

For $i \in \{1, \ldots, m\}$, let $A_i$, $A^-_i$, $A^+_i$ denote the $i^{th}$ row of A, $A^-$, $A^+$ respectively.

Clearly, $A^-_i \leq A^+_i$ for all $i \in \{1, \ldots, m\}$.

Suppose that for all $i \in \{1, \ldots, m\}$ it is the case that: (i) $A^-_i \leq A_i \leq A^+_i$, and (ii) there exists $x^* \geq 0$ such that $A_i x^* = b_i$.

If $x^*$ solves the LP problem [Maximize $p^T x$, subject to $A^-_i x \leq b_i \leq A^+_i x$ for all $i \in \{1, \ldots, m\}$, $x \geq 0$], then $x^*$ solves P1.

This follows from the observation that for all $x \geq 0$ and $i \in \{1, \ldots, m\}$: if $A_i x = b_i$, then $A^-_i x \leq b_i \leq A^+_i x$.

However, the converse - in general- may not be true.

Consider the linear system $Ax = b$, $x \geq 0$.

Let X be the set of all solutions of this linear system.

In what follows we assume that $X \neq \phi$.

Clearly the set X is closed and convex.

Let $X^*$ be the set $\{x \in X|$ the columns of A in the array $\langle A^j| x_j > 0 \rangle$ are linearly independent$\}$.

Elements of $X^*$ are called **basic solutions**.

The relationship between X and $X^*$ will be our primary concern in section 6 of this paper.

**3. A Non-Substitution Theorem for Linear Programming:** In this section we present a very general sensitivity result.

For what follows in this section, suppose, $\bar{x}$ is <u>a solution of problem</u> P1.

**Theorem 1:** If an n-dimension vector $x^* \geq 0$ satisfies $\{j|x_j^* > 0\} \subset \{j|\bar{x}_j > 0\}$, then $x^*$ solves the following LP problem, which we shall refer to as P2.

Maximize $p^T x$, subject to $Ax = Ax^*$, $x \geq 0$.

**Proof:** Let $b^*$ denote $Ax^*$.

The dual of P1 is defined as follows:

Minimize $y^T b$ subject to $y^T A \geq p^T$, y unconstrained in sign.

By the Weak duality theorem, the dual has a solution $\bar{y}$ and $\bar{y}^T b = p^T \bar{x}$.

By the Complementary Slackness Theorem, $\{j| \bar{x}_j > 0\} \subset \{j|\bar{y}^T A^j = p_j\}$.

The dual of P2 is defined as follows:

Minimize $y^T b^*$ subject to $y^T A \geq p^T$, y unconstrained in sign.

Clearly $\bar{y}$ satisfies the constraints of this problem.

Further, $(\bar{y}^T A - p^T)x^* = \sum_{j=1}^n (\bar{y}^T A^j - p_j)x_j^* = \sum_{\{j| x_j^* > 0\}}(\bar{y}^T A^j - p_j)x_j^* + \sum_{\{j|x_j^* = 0\}}(\bar{y}^T A^j - p_j)x_j^* = \sum_{\{j| x_j^* > 0\}}(\bar{y}^T A^j - p_j)x_j^*$.

However, $\{j| x_j^* > 0\} \subset \{j| \bar{x}_j > 0\} \subset \{j|(\bar{y}^T A^j = p_j)\}$.

Thus, $\sum_{\{j| x_j^* > 0\}}(\bar{y}^T A^j - p_j)x_j^* = 0$ and hence $(\bar{y}^T A - p^T)x^* = 0$.

However, $Ax^* = b^*$ and $(\bar{y}^T A - p^T)x^* = 0$ implies $\bar{y}^T b^* = p^T x^*$.

Since $x^*$ satisfies the constraints of P2, $\bar{y}$ satisfies the constraints for its dual and $\bar{y}^T b^* = p^T x^*$, it must be the case that $x^*$ solves P2 and $\bar{y}$ solves its dual. Q.E.D.

**4. A proposition followed by a proof of theorem 1 in Mangasarian (1979):** For what follows in this section, suppose, $\bar{x}$ is <u>a solution of problem</u> P1.

We present below a lemma with a proof that is easy to follow. However, before that we need to define a few concepts.

Let $I^0(\bar{x}) = \{j|\bar{x}_j = 0\}$. This set will play a very important role in what follows. We will refer to $I^0(\bar{x})$ as **the set of zero-coordinates of $\bar{x}$**.

**Notation:** For all n-dimensional column vectors h, let $\|h\| = \sqrt{\sum_{j=1}^{n} h_j^2}$, where for all $j \in \{1, \ldots, n\}$, $h_j$ is the $j^{th}$ coordinate of h.

**Lemma 1:** The problem: Maximize $p^T h$, subject to $Ah = 0$, $h_j \geq 0$ for all $j \in I^0(\bar{x})$, has the 0-vector as its unique solution

<u>if and only if</u>

for all $q \neq 0$, there exists $\varepsilon > 0$ such that for all $\delta \in (0, \varepsilon)$, the 0-vector, is <u>a</u> solution for the problem: Maximize $(p + \delta q)^T h$, subject to $Ah = 0$, $h_j \geq 0$ for all $j \in I^0(\bar{x})$.

**Proof:** Suppose that the problem: Maximize $p^T h$, subject to $Ah = 0$, $h_j \geq 0$ for all $j \in I^0(\bar{x})$, has the 0-vector as its unique solution.

Towards a contradiction suppose that for some $q \neq 0$, it is the case that for all $\varepsilon > 0$, there exists $\delta \in (0, \varepsilon)$ such that $(p + \delta q)^T h > 0$, $Ah = 0$, $h_j \geq 0$ for all $j \in I^0(\bar{x})$, has a non-zero solution.

Thus, for all $\varepsilon > 0$, there exists $\delta \in (0, \varepsilon)$ such that $(p + \delta q)^T h > 0$, $Ah = 0$, $h_j \geq 0$ for all $j \in I^0(\bar{x})$ $\|h\| = 1$, has a solution.

Thus, there is a strictly decreasing sequence $\langle \delta(k) | k \in \mathbb{N} \rangle$ such that $\lim_{k \to \infty} \delta(k) = 0$ and a sequence $\langle h(k) | k \in \mathbb{N} \rangle$ in $\mathbb{R}_+^n \setminus \{0\}$, such that $(p + \delta(k)q)^T h(k) > 0$, $Ah(k) = 0$, $h_j(k) \geq 0$ for all $j \in I^0(\bar{x})$, $\|h(k)\| = 1$ for all $k \in \mathbb{N}$.

$\langle h(k) | k \in \mathbb{N} \rangle$ is a bounded infinite sequence of real numbers, since it lies on the surface of the unit sphere.

The Bolzano-Weierstrass theorem, which is a well-known result in basic real analysis (Proposition 6, in Topic 5 of Lahiri (2020)), says that $\langle h(k) | k \in \mathbb{N} \rangle$ has a convergent "subsequence" and since $Ah(k) = 0$, $h_j(k) \geq 0$ for all $j \in I^0(\bar{x})$, $\|h(k)\| = 1$ for all $k \in \mathbb{N}$, the limit, say $\bar{h}$, that the subsequence converges to must satisfy $A\bar{h} = 0$, $\bar{h}_j \geq 0$ for all $j \in I^0(\bar{x})$ and $\|\bar{h}\| = 1$.

Since $\lim_{k \to \infty} \delta(k) = 0$ and since $(p + \delta(k)q)^T h(k) > 0$ for all $k \in \mathbb{N}$, the subsequential limit of the sequence $\langle (p + \delta(k)q)^T h(k) | k \in \mathbb{N} \rangle$ which is $p^T \bar{h}$ is non-negative, i.e., $p^T \bar{h} \geq 0$.

Since $\|\bar{h}\| = 1$, the problem: Maximize $p^T h$, subject to $Ah = 0$, $h_j \geq 0$ for all $j \in I^0(\bar{x})$, has a non-zero solution, contradicting our assumption that '0' is the *unique* solution for this problem.

Now. suppose, $\bar{h} \neq 0$ is a solution for the problem: Maximize $p^T h$, subject to $Ah = 0$, $h_j \geq 0$ for all $j \in I^0(\bar{x})$.

Thus, $A\bar{h} = 0$, $\bar{h}_j \geq 0$ for all $j \in I^0(\bar{x})$ and $p^T \bar{h} \geq 0$.

Since $\bar{h} \neq 0$, there exists $\bar{q} \in \mathbb{R}^n \setminus \{0\}$ such that $\bar{q}^T \bar{h} > 0$. In fact, $\bar{q} = \bar{h}$, satisfies $\bar{q}^T \bar{h} = \bar{h}^T \bar{h} > 0$, since $\bar{h} \neq 0$.

Thus, $(p + \alpha \bar{q})^T \bar{h} > 0$ for all $\alpha > 0$.

This, violates the requirement that for all non-zero n-dimensional column vectors q, there exists $\varepsilon > 0$ such that for all $\delta \in (0, \varepsilon)$, the 0-vector, is a solution for the second maximization problem in the statement of the lemma.

This proves the lemma. Q.E.D.

An almost immediate consequence of lemma 1 is the following result.

**Proposition 1:** The following three statements are equivalent:

(i) $\bar{x}$ is the <u>unique</u> solution for P1, i.e., Maximize $p^T x$, subject to $Ax = b$, $x \geq 0$.

(ii) The 0-vector is the unique solution of the problem: Maximize $p^T h$, subject to $Ah = 0$, $h_j \geq 0$ for all $j \in I^0(\bar{x})$.

(iii) For all $q \neq 0$, there exists $\varepsilon > 0$ such that for all $\delta \in (0, \varepsilon)$, the 0-vector, is <u>a</u> solution for the problem: Maximize $(p + \delta q)^T h$, subject to $Ah = 0$, $h_j \geq 0$ for all $j \in I^0(\bar{x})$.

**Proof:** Since lemma 1 establishes the equivalence of (ii) and (iii) it is enough to establish the equivalence of (i) and (ii).

This is true because:

$[\bar{h} \neq 0, A\bar{h} = 0, \bar{h}_j \geq 0$ for all $j \in I^0(\bar{x}), p^T \bar{h} \geq 0]$ 'if and only if' $[\bar{x} + \varepsilon \bar{h} \neq \bar{x}, A(\bar{x} + \varepsilon \bar{h}) = b$, $p^T(\bar{x} + \varepsilon \bar{h}) \geq p^T \bar{x}$ for all $\varepsilon > 0$ sufficiently small satisfying $\bar{x}_j + \varepsilon \bar{h}_j \geq 0$ for all $j \in \{1,\ldots, n\}]$.

This establishes the equivalence of (i) and (ii). Q.E.D.

Theorem 1 of Mangasarian (1979) now follows as an immediate consequence of Proposition 1.

**Theorem 2 (Theorem 1 of Mangasarian (1979)):** $\bar{x}$ is the <u>unique</u> solution of P1, i.e., the problem: Maximize $p^T x$, subject to $Ax = b$, $x \geq 0$

<u>if and only if</u>

for all $q \neq 0$, there exists $\varepsilon > 0$ such that for all $\delta \in (0, \varepsilon)$, $\bar{x}$ is 'a' solution for:

Maximize $(p + \delta q)^T x$, subject to $Ax = b$, $x \geq 0$.

**Proof:** In view of Proposition 1 all that we need to show is the following result.

For all $q \neq 0$, there exists $\varepsilon > 0$ such that for all $\delta \in (0, \varepsilon)$, the following two statements are equivalent:

(a) $h = 0$, is <u>a</u> solution of: Maximize $(p + \delta q)^T h$, subject to $Ah = 0$, $h_j \geq 0$ for all $j \in I^0(\bar{x})$.

(b) $\bar{x}$ is <u>a</u> solution of: Maximize $(p + \delta q)^T x$, subject to $Ax = b$, $x \geq 0$.

Suppose statement (a) holds. Towards a contradiction suppose that there exists $x \geq 0$ such that $Ax = b$ and $(p + \delta q)^T x > (p + \delta q)^T \bar{x}$.

Let $h = x - \bar{x} \neq 0$. Since $x_j \geq 0 = \bar{x}_j$ for all $j \in I^0(\bar{x})$, it must be the case that $h_j \geq 0$, for all $j \in I^0(\bar{x})$,

Thus, $Ah = 0$ and $(p + \delta q)^T h > 0$, $h_j \geq 0$, for all $j \in I^0(\bar{x})$ contradicts the assumption that the 0-vector solves: Maximize $(p + \delta q)^T h$, subject to $Ah = 0$, $h_j \geq 0$ for all $j \in I^0(\bar{x})$.

Thus, statement (b) must hold.

Now suppose that statement (b) holds. Towards a contradiction suppose there exists h satisfying $Ah = 0$, $h_j \geq 0$, for all $j \in I^0(\bar{x})$ and $(p + \delta q)^T h > 0$.

For $\varepsilon > 0$ sufficiently small, $\bar{x} + \varepsilon h \geq 0$ and $A(\bar{x} + \varepsilon h) = A\bar{x} + \varepsilon Ah = b$. Further, $(p + \delta q)^T (\bar{x} + h) = (p + \delta q)^T \bar{x} + (p + \delta q)^T h > (p + \delta q)^T \bar{x}$, contradicting that (b) holds.

Thus, statement (a) must hold. Q.E.D.

**5. Some observations:** It is worth noting that section 4 of Szilagyi (2006), is also concerned with a given solution for a LP problem being its unique solution. It appears, that in Szilagyi (2006), the given solution is assumed to be a "non-degenerate" basic feasible solution. The primal problem that the paper is concerned with has constraints that are inequalities, instead of equalities, the latter being the way we have formulated our problem. Further, unlike here where the primal variables are assumed to be non-negative, the primal variables there are unconstrained in sign. The two theorems 4.1 and 4.3 in Szilagyi (2006)- these two theorems being the "work horse" for the other results in section 4 of Szilagyi (2006)- are analogous to the equivalence of the first two statements in our proposition 1 in this paper. Our theorem 2 is, as already stated, available in Mangasarian (1979), and our proposition 1 helps us to provide an alternative (and hopefully new?) proof of the same.

In a 2017 paper by Appa (Appa (2017)), which is devoted entirely to the uniqueness problem, there is a necessary and sufficient condition for an optimal solution for a LP problem to be its unique solution. Unlike our analysis here, the "candidate" optimal solution in Appa (2017) is assumed to be a basic optimal solution. While, basic optimality is implied by the uniqueness of the candidate optimal solution, the converse is not true. Lemma 1 in Appa (2017) says: given $\bar{x}$ is a solution for P1 – not necessarily basic- an n-dimensional vector $\tilde{x}$ is a solution for P1 <u>if and only if</u> $\tilde{x}$ is a solution for the linear system, $Ax = b$, $p^T x = p^T \bar{x}$, $x \geq 0$. The result concerning uniqueness in Appa (2017) is an algorithm that requires solving a LP problem whose constraints are the same as in Lemma 1 in Appa (2017), but the coefficients of the objective function are defined as follows: the coefficients corresponding to the basic variables in the candidate solution are '0' and the coefficients corresponding to the non-basic variables in the candidate solution are '1'. The main result in Appa (2017) states, that the candidate optimal solution is unique <u>if and only if</u> the optimal value of this LP problem is '0'. Thus, if $\bar{x}$ is a solution for P1 and in addition $\bar{x}$ is basic, then $\bar{x}$ is "the unique solution" for P1 <u>if and only if</u> the n-dimensional '0' is the optimal value of the following LP problem: Maximize $q^T x$, subject to $Ax = b$, $p^T x = p^T \bar{x}$, $x \geq 0$, where q is an-dimensional vector satisfying $q^T \bar{x} = 0$ and $q_j = 1$ for all j such that $\bar{x}_j = 0$.

Our results so far, are in no way related to any result or its implementation in Appa (2017).

This main result in Appa (2017) and our proposition 1, suggests an alternative necessary and sufficient condition for uniqueness of a given optimal solution.

**Alternative to Appa's (2017) theorem:** If $\bar{x}$ is an optimal solution for P1 and in addition $\bar{x}$ is basic, then $\bar{x}$ is the unique optimal solution of P1 <u>if and only if</u> '0' is the optimal value of the following linear programming problem: Maximize $q^T h$, subject to $Ah = 0$, $p^T h = 0$, $h \geq 0$ for all $j \in I^0(\bar{x})$, where q is an-dimensional vector satisfying $q^T \bar{x} = 0$ and $q_j = 1$ for all j such that $\bar{x}_j = 0$.

The proof of the alternative that we suggest, is quite easy.

**6. Convex combinations of basic feasible solutions:** Recall from section 2, that given a system of linear equations $Ax = b$, $x \geq 0$, where A is an m×n matrix of real numbers and b is an m-dimensional column vector, we use X to denote the set of all solutions of this linear system.

In what follows we assume that $X \neq \phi$. Clearly the set X is closed and convex.

Recall the definition of basic solutions $X^*$ in section 2: $X^* = \{x \in X | <A^j | x_j > 0>$ are linearly independent$\}$.

By theorem 6 in Topic 1 of Lahiri (2020) called the "Non-negative Pivot Theorem", we get that $X^* \neq \phi$.

Since the collection of sets of linearly independent columns of A is finite, $X^*$ must be finite.

If $0 \in X$, then it is a basic solution and hence $0 \in X^*$.

A point $x \in X$, is said to be an **extreme point** of X if there does not exist any pair of distinct n-vectors $x^{(1)}, x^{(2)} \in X$, such that $x = \alpha x^{(1)} + (1-\alpha) x^{(2)}$ for some $\alpha \in (0, 1)$.

**Lemma 2:** Let $x \in X^*$. Then x is an extreme point of X.

**Proof:** Let $x \in X^*$ and towards a contradiction suppose there exists a pair of distinct n-vectors $x^{(1)}, x^{(2)} \in X$, such that $x = \alpha x^{(1)} + (1-\alpha) x^{(2)}$ for some $\alpha \in (0, 1)$.

Thus, $x^{(1)} \neq x^{(2)}$.

Since $x^{(1)}, x^{(2)} \geq 0$ and x is a convex combination of $x^{(1)}, x^{(2)}$, it must be the case that $x_j = 0$ implies $x_j^{(k)} = 0$ for all $k = 1, 2$.

Thus, $\sum_{\{j | x_j > 0\}} A^j x_j^{(k)} = b$ for all $k = 1, 2$.

Since $x^{(1)} \neq x^{(2)}$, $\sum_{\{j | x_j > 0\}} A^j (x_j^{(1)} - x_j^{(2)}) = 0$ with $x_j^{(1)} - x_j^{(2)} \neq 0$ for at least one j satisfying $x_j > 0$.

Thus, the columns in $<A^j | x_j > 0>$ are linearly dependent, contradicting $x \in X^*$. Q.E.D.

Recall the linear optimization problem with linear objective function P1:

[Maximize $p^T x$, subject to $Ax = b$, $x \geq 0$] where p is some n-dimensional vector.

Suppose the *set of optimal solutions* Y for this problem is non-empty. It is easily verified that Y is a closed and convex subset of $\mathbb{R}_+^n$

Let $Y^* = Y \cap X^*$. Points in $Y^*$ are referred to as **basic optimal solution**.

Note that $Y^* \subset X^*$.

By Proposition 4 in Topic 2 of Lahiri (2020), $Y^*$ is non-empty and finite.

A point $y \in Y$, is said to be an **extreme point** of Y if there does not exist any pair of distinct n-vectors $y^{(1)}, y^{(2)} \in Y$, such that $y = \alpha y^{(1)} + (1-\alpha) y^{(2)}$ for some $\alpha \in (0, 1)$.

The application of the results in R4.4 of Lancaster (1968) in "linear programming" that we referred to in the introductory section of this paper, is very likely the following.

**Theorem 3:** Suppose X is bounded. Then $Y = \{\sum_{y^* \in Y^*} \alpha_{y^*} y^* | \alpha_{y^*} \geq 0$ for all $y^* \in Y^*$ and $\sum_{y^* \in Y^*} \alpha_{y^*} = 1\}$ with $Y^*$ as the set of extreme points of Y.

**Proof:** It is easy to see that $\sum_{y^* \in Y^*} \alpha_{y^*} y^* | \alpha_{y^*} \geq 0$ for all $y^* \in Y^*$ and $\sum_{y^* \in Y^*} \alpha_{y^*} = 1\} \subset Y$.

Suppose X is bounded and towards a contradiction suppose $Y \setminus \{\sum_{y^* \in Y^*} \alpha_{y^*} y^* | \alpha_{y^*} \geq 0$ for all $y^* \in Y^*$ and $\sum_{y^* \in Y^*} \alpha_{y^*} = 1\} \neq \phi$. Let $y \in Y \setminus \{\sum_{y^* \in Y^*} \alpha_{y^*} y^* | \alpha_{y^*} \geq 0$ for all $y^* \in Y^*$ and $\sum_{y^* \in Y^*} \alpha_{y^*} = 1\}$.

For all $y^* \in Y^*$ and $x^* \in X^* \setminus Y^*$, it must be that $p^T y^* > p^T x^*$.

Since $p^T y > p^T x^*$ for all $x^* \in X^* \setminus Y^*$: $[y = \sum_{x^* \in X^*} \alpha_{x^*} x^*, \alpha_{x^*} \geq 0$ for all $x^* \in X^*$, $\alpha_{x^*} > 0$ for some $x^* \in X^* \setminus Y^*$ and $\sum_{x^* \in X^*} \alpha_{x^*} = 1]$ implies $[p^T y = \sum_{x^* \in Y^*} \alpha_{x^*} p^T x^* < p^T y^*$ for all $y^* \in Y^*]$, contradicting $y \in Y$.

Thus, $[y \in Y \setminus \{\sum_{y^* \in Y^*} \alpha_{y^*} y^* | \alpha_{y^*} \geq 0$ for all $y^* \in Y^*$ and $\sum_{y^* \in Y^*} \alpha_{y^*} = 1\}]$ implies that the linear system $[y = \sum_{x^* \in X^*} \alpha_{x^*} x^*, \alpha_{x^*} \geq 0$ for all $x^* \in X^*$ and $\sum_{x^* \in X^*} \alpha_{x^*} = 1]$ has no solution.

By Farkas' lemma, there exists an n-vector q and a real number $\beta$ such that $q^T x^* + \beta \leq 0$ for all $x^* \in X^*$ and $q^T y + \beta > 0$. Thus, $q^T y > q^T x^*$ for all $x^* \in X^*$. Clearly for such a strict inequality to hold it is necessary that $q \neq 0$.

Since X is bounded, by Proposition 4, in Topics 2 of Lahiri (2024)), it follows that the linear optimization problem [Maximize $q^T x$, subject to $Ax = b$, $x \geq 0$] has at least one (optimal) solution that is basic, i.e., it has at least one optimal solution that belongs to $X^*$.

However, $y \in Y \subset X$ and $q^T y > q^T x^*$ for all $x^* \in X^*$ implies that $(\underset{x \in X}{\arg\max}\, q^T x) \cap X^* = \phi$, leading to a contradiction.

Hence, it must be the case that $Y \setminus \{\sum_{y^* \in Y^*} \alpha_{y^*} y^* | \alpha_{y^*} \geq 0$ for all $y^* \in Y^*$ and $\sum_{y^* \in Y^*} \alpha_{y^*} = 1\} = \phi$.

Since $\{\sum_{y^* \in Y^*} \alpha_{y^*} y^* | \alpha_{y^*} \geq 0$ for all $y^* \in Y^*$ and $\sum_{y^* \in Y^*} \alpha_{y^*} = 1\} \subset Y$, it follows that $Y = \{\sum_{y^* \in Y^*} \alpha_{y^*} y^* | \alpha_{y^*} \geq 0$ for all $y^* \in Y^*$ and $\sum_{y^* \in Y^*} \alpha_{y^*} = 1\}$.

Since $Y^* = Y \cap X^*$, it follows from lemma 1 that $Y^*$ is the set of extreme points of Y. Q.E.D.

**Note:** One reason for the relative scarcity of elementary (or simple) proofs of results in linear optimization such as theorem 3, may be that except for Mote and Madhavan (2016), most contemporary work in this area is either "obsessively algorithmic" or relies entirely on the "simplex method" (to the total exclusion of all else) as its "work horse".

It is important to note that the assumption on X being bounded is required as the following example reveals.

**Example:** Consider the problem [Maximize $y_1$, subject to $[1\ 0] \begin{bmatrix} y_1 \\ y_2 \end{bmatrix} = 1$, $y_1, y_2 \geq 0$].

In this case $X = \{(1, y_2) | y_2 \geq 0\}$ which is unbounded and $X^* = \{\begin{bmatrix} 1 \\ 0 \end{bmatrix}\}$.

Clearly, $Y = X$ although $Y^* = Y \cap X^* = X^* = \{\begin{bmatrix} 1 \\ 0 \end{bmatrix}\}$

Further, $\begin{bmatrix} 1 \\ 1 \end{bmatrix} \in Y$ and it cannot be obtained as a convex combination of solutions in $Y^*$.

By taking p = 0 in the definition of the linear optimization problem with linear objective function that we are concerned with in this note, we get the following result.

**Corollary of Theorem 3:** Suppose X is bounded. Then $X = \{\sum_{x^* \in X^*} \alpha_{x^*} x^* | \alpha_{x^*} \geq 0$ for all $x^* \in X^*$ and $\sum_{x^* \in X^*} \alpha_{x^*} = 1\}$ with $X^*$ as the set of extreme points of X.

**Proof:** Follows immediately from the fact that if p = 0, then $Y = X$ and $Y^* = X^*$. Q.E.D.

**7. Proof of "The Birkhoff- von Neumann Theorem":** We provide here a proof of "The Birkhoff- von Neumann Theorem" as an application of Corollary of Theorem 3, and the proposition that every doubly stochastic matrix is a convex combination of permutation matrices.

A "non-negative matrix" P is said to be **doubly stochastic** if all row sums of P are 1 and all column sums of P are 1.

Thus, the sum of all column sums of P = the number of rows of P, and the sum of all row sums of P = the number of columns of P.

However, the sum of all column sums of P = the sum of all entries in the matrix P = the sum of all row sums of P.

Thus, a doubly stochastic matrix must have the same number of rows and columns, i.e., it must be a square matrix.

A doubly stochastic matrix all whose entries are 0 or 1, is said to be a **permutation matrix**.

Thus, each row of a permutation matrix has exactly one entry that is 1 and all other entries in the row are 0, and each column of a permutation matrix has exactly one entry that is 1 and all other entries in the column are 0.

Thus, permutation matrices are obtained by permuting the columns of the identity matrix and hence given any positive integer 'n', the total number of n×n permutation matrices is n!.

Let $\{P^{(k)} | k \in \{1, \ldots, n!\}\}$ be the set of n×n permutation matrices.

Part (i) of the following result has a very elegant proof that is available in section 2 of chapter 8 in Hurlbert (2010) which for the sake of being self-contained is being reproduced below.

**Proposition 2:** (i) If a n×n doubly stochastic matrix P is not a permutation matrix, then there exists two distinct (i.e., unequal) doubly stochastic matrices $Q^{(1)}$, $Q^{(2)}$ such that $P = \frac{1}{2}Q^{(1)} + \frac{1}{2}Q^{(2)}$.

(ii) If P is a permutation matrix, then any convex combination of matrices in $\{P^{(k)}| k \in \{1, \ldots, n!\}\}$, must assign weight 1 to P and weight 0 to all other matrices in P.

**Proof:** (i) Suppose P is a n×n doubly stochastic matrix that is not a permutation matrix. Then then there exists an entry in P, say the entry at the intersection row $r_1$ and column $s_1$ denoted by $p(r_1, s_1)$, such that $0 < p(r_1, s_1) < 1$.

Since P is a doubly stochastic matrix, the sum of entries in each row is 1, and so there exists another entry in row $r_1$, say $p(r_1, s_2)$ with $s_2 \neq s_1$ such that $0 < p(r_1, s_2) < 1$.

Since the sum of entries in each column is 1, there exists another entry in column $s_2$, say $p(r_2, s_2)$ with $r_2 \neq r_1$ such that $0 < p(r_2, s_2) < 1$.

Having found $\{p(r_i, s_i)| i = 1, \ldots, j\} \cup \{p(r_j, s_{i+1})| i = 1, \ldots, j\}$ for some $j \geq 1$, all distinct, since the sum of entries in column j+1 is equal to 1, there exists a row $r_{j+1}$, such that $0 < p(r_{j+1}, s_{j+1}) < 1$.

Similarly, having found $\{p(r_i, s_i)| i = 1, \ldots, j\} \cup \{p(r_{i-1}, s_i)| i = 2, \ldots, j\}$ for some $j \geq 2$, all distinct, since the sum of entries in row j is equal to 1, there exists a column $s_{j+1}$, such that $0 < p(r_j, s_{j+1}) < 1$.

Since, the matrix P has a finite number of entries, this process cannot go on for ever without revisiting a row-column pair that was visited earlier in the process, and so there is a first time when the process returns to a row-column pair which it had either started out from or visited earlier, thereby forming a cycle.

Consider such a cycle of minimal length beginning with some $p(r_1, s_1)$ satisfying $0 < p(r_1, s_1) < 1$ and terminating at row $r_k = r_1$.

If $(r_k, s_{k+1}) = (r_1, s_1)$, then $(r_1, s_k) = (r_k, s_k) \neq (r_k, s_{k+1}) = (r_1, s_1) \neq (r_1, s_2)$, and so the cycle could be decreased in length by removing $(r_1, s_1)$ and going directly from $(r_1, s_k)$ to $(r_1, s_2)$, thereby contradicting it has minimal length.

Thus, a cycle of minimal length must satisfy $(r_k, s_k) = (r_1, s_1)$, i.e. includes an even number of distinct row-column pairs.

Let, $\varepsilon_0 = \min \{p(r_j, s_j), 1 - p(r_j, s_j), p(r_j, s_{j+1}), 1 - p(r_j, s_{j+1})| j = 1, \ldots, k\}\}$, where $s_{k+1} = s_1$. Clearly $0 < \varepsilon_0 < 1$.

Let $\varepsilon$ be a real number satisfying $0 < \varepsilon < \varepsilon_0$.

Let $Q^{(1)}$ be the n×n matrix such that $q^{(1)}(r, s) = p(r,s)$ if $(r,s) \notin \{(r_j, s_j)| j = 1, \ldots, k\} \cup \{(r_{j-1}, s_j)| j = 2, \ldots, k\}$, $q^{(1)}(r, s) = p(r,s) - \varepsilon$ if $(r,s) \in \{(r_j, s_j)| j = 1, \ldots, k\}$ and $q^{(1)}(r, s) = p(r,s) + \varepsilon$ if $(r,s) \in \{(r_{j-1}, s_j)| j = 2, \ldots, k\}$.

Let $Q^{(2)}$ be the n×n matrix such that $q^{(2)}(r, s) = p(r,s)$ if $(r,s) \notin \{(r_j, s_j)| j = 1, \ldots, k\} \cup \{(r_{j-1}, s_j)| j = 2, \ldots, k\}$, $q^{(2)}(r, s) = p(r,s) + \varepsilon$ if $(r,s) \in \{(r_j, s_j)| j = 1, \ldots, k\}$ and $q^{(2)}(r, s) = p(r,s) - \varepsilon$ if $(r,s) \in \{(r_{j-1}, s_j)| j = 2, \ldots, k\}$.

Thus, $Q^{(1)} \neq Q^{(2)}$ and both are doubly stochastic. Further, $P = \frac{1}{2}Q^{(1)} + \frac{1}{2}Q^{(2)}$.

(ii) Suppose P is a permutation matrix and suppose $P = \sum_{k=1}^{n!} \alpha_k P^{(k)}$ with $\alpha_k \geq 0$ for all $k = 1, \ldots, n!$ and $\sum_{k=1}^{n!} \alpha_k = 1$.

Suppose that the column in row 1 of P that has 1, is $s_1$. There are (n-1)! permutation matrices that have 1 in column $s_1$ of row 1. Let $S_1$ be the set of such matrices. Clearly, $\alpha_k = 0$ for all permutation matrices that do not belong to $S_1$ since P being a permutation matrix can have at most one entry that is 1, in each row. Now suppose that the column in row 2 of P that has 1, is $s_2$. There are (n-2)! permutation matrices in $S_1$ that have 1 in column $s_2$ of row 2. Let $S_2$ be the subset of such matrices in $S_2$. Clearly, $\alpha_k = 0$ for all permutation matrices that do not belong to $S_2$. Proceeding in this manner we arrive at the set $S_{n-1}$ of permutation matrices that in each of the first n-1 rows has 1 in the same column as P has, and the cardinality of $S_{n-1}$ is (n- (n-1))! = 1! = 1. Further, $\alpha_k = 0$ for all permutation matrices that do not belong to $S_{n-1}$. Thus, the first n-1 rows of the permutation matrix in $S_{n-1}$ are equal to the first n-1 rows of P. Since both P and the one in $S_{n-1}$ are permutation matrices, their first n-1 rows uniquely determine their last row. Thus $S_{n-1} = \{P\}$ and $\alpha_k = 1$ for the unique k for which $P = P^{(k)}$.
Q.E.D.

We now introduce some notations.

For any n×n matrix B let $\mathcal{A}(B)$ be the $n^2$ dimension column vector such the for each $j \in \{1, \ldots, n\}$, the coordinates of $\mathcal{A}(B)$ numbered $(j-1)n + 1, \ldots, jn$ form the column vector $B^j$ which is the $j^{th}$ column of B.

For the positive integer n, let $E^{(n)}$ denote the n-dimensional column sum vector, i.e., the n-dimensional column vector all whose coordinates are 1 and for all $i \in \{1, \ldots, n\}$, let $E^{(n,i)}$ denote the n-dimensional column unit vector, i.e., the n-dimensional column vector whose $i^{th}$ coordinate is 1 and all other coordinates are 0.

Let $\mathcal{E}^{(1)}$ be the n×n² matrix, such that for each $j \in \{1, \ldots, n\}$, the coordinates in the $j^{th}$ row of $\mathcal{E}^{(1)}$ occupying columns $(j-1)n + 1, \ldots, jn$ is the row vector $E^{(n)T}$ and all other entries in row j are 0.

Let $\mathcal{E}^{(2)}$ be the n×n² matrix, such that for each $j \in \{1, \ldots, n\}$, the coordinates in the $j^{th}$ row of $\mathcal{E}^{(1)}$ occupying columns $(h-1)n + 1, \ldots, hn$ is the row vector $E^{(n,j)T}$ for all $h \in \{1, \ldots, n\}$.

Thus, a n× n matrix B is doubly stochastic <u>if and only if</u> $\mathcal{A}(B) \in \mathbb{R}_+^{n^2}$ and $\begin{bmatrix} \mathcal{E}^{(1)} \\ \mathcal{E}^{(2)} \end{bmatrix} \mathcal{A}(B) = \begin{bmatrix} E^{(n)} \\ E^{(n)} \end{bmatrix}$.

Let X be the set of all n × n doubly stochastic matrices, i.e., $X = \{\mathcal{A} \in \mathbb{R}_+^{n^2} | \begin{bmatrix} \mathcal{E}^{(1)} \\ \mathcal{E}^{(2)} \end{bmatrix} \mathcal{A} = \begin{bmatrix} E^{(n)} \\ E^{(n)} \end{bmatrix}\}$.

A pictorial discussion of the following theorem is available here:
https://www.theoremoftheday.org/CombinatorialTheory/Birkhoff/TotDBirkhoff.pdf

**Theorem 4 (Birkhoff (1946) - von Neumann (1953)):** A n×n matrix is doubly stochastic <u>if and only if</u> it can be expressed as a convex combination of n×n permutation matrices.

**Proof:** From Proposition 2, we know that $X^* \subset \{\mathcal{A}(P^{(k)})| k = 1, …, n!\} \subset X$.

Further, it is easily verified that X is a bounded set.

Thus, from Corollary of theorem 3, it follows that $X = \{\sum_{\mathcal{A}^* \in X^*} \alpha_{\mathcal{A}^*}\mathcal{A}^*| \alpha_{\mathcal{A}^*} \geq 0$ for all $\mathcal{A}^* \in X^*$ and $\sum_{\mathcal{A}^* \in X^*} \alpha_{\mathcal{A}^*} = 1\}$.

However, $\{\sum_{\mathcal{A}^* \in X^*} \alpha_{\mathcal{A}^*}\mathcal{A}^*| \alpha_{\mathcal{A}^*} \geq 0$ for all $\mathcal{A}^* \in X^*$ and $\sum_{\mathcal{A}^* \in X^*} \alpha_{\mathcal{A}^*} = 1\} \subset \{\sum_{k=1}^{n!} \alpha_k \mathcal{A}(P^{(k)})| \alpha_k \geq 0$ for all $k = 1, …, n!$ and $\sum_{k=1}^{n!} \alpha_k = 1\} \subset X$.

Thus, $X = \{\sum_{k=1}^{n!} \alpha_k \mathcal{A}(P^{(k)})| \alpha_k \geq 0$ for all $k = 1, …, n!$ and $\sum_{k=1}^{n!} \alpha_k = 1\}$. Q.E.D.

An immediate consequence of theorem 4 and part (ii) of proposition 2 is the following corollary.

**Corollary of Theorem 4:** $X^* = \{\mathcal{A}(P^{(k)})| k = 1, …, n!\}$, i.e., the set of extreme points of the set of all n×n doubly stochastic matrices is the set of all n×n permutation matrices.

A different application of LP is used in a proof of theorem 4, that is available in Azar (2010).

**Acknowledgment:** This paper is a revised extension of an earlier version entitled "The Non-Substitution Theorem and Uniqueness of Solution in Linear Programming". I wish to thank Paolo Serafini and Jose A. Silva for comments and clarifications about results cited and discussed in this paper. I also wish to thank Glenn Hurlbert for providing me with the updated bibliographic information about his work that has been cited here and for suggesting an improvement in the proof of part (i) of Proposition 2 in this paper. I want to put on record my grateful appreciation of Sci-Hub (https://sci-hub.se/) for providing me on October 5, 2024, with free access to Appa (2017) that is cited in this paper.